\documentclass[a4paper, 12pt]{article}

\usepackage{amssymb,amsbsy,amsmath,amsfonts,amssymb,amscd}
\usepackage{latexsym,eucal,exscale,epic,eepic,epsfig,xcolor}
\usepackage{mathrsfs}

\newtheorem{theorem}{Theorem}[section]
\newtheorem{lemma}[theorem]{Lemma}
\newtheorem{proposition}[theorem]{Proposition}
\newtheorem{corollary}[theorem]{Corollary}

\newtheorem{conjecture}[theorem]{Conjecture}

\newcommand{\ds}{\displaystyle}
\newcommand{\GG}{\mathbb G}
\newcommand{\ofr}{{\mathfrak o}}
\newcommand{\pfr}{{\mathfrak p}}

\newcommand{\lra}{\longrightarrow}
\newcommand{\noi}{\noindent}
\newcommand{\Ap}{\mathscr A}

\newcommand{\St}{{\mathbf S}{\mathbf t}}

\newcommand{\ZZ}{\mathbb Z}

\newcommand{\Sm}{\mathscr S}
\newcommand{\CC}{\mathbb C}

\newcommand{\HH}{{\mathbb H}}
\newcommand{\HIw}{{\mathscr H}}
\newcommand{\Iw}{{\rm I}}
\newcommand{\VV}{\mathscr V}
\newcommand{\FF}{\mathscr F}
\newcommand{\bV}{{\mathbf V}}

\newcommand{\Ch}{{\mathscr C \! h}}

\newcommand{\MM}{\mathbf M}
\newcommand{\bM}{\mathbf M}

\newcommand{\WAff}{W^{\rm Aff}}
\newcommand{\SL}{{\rm SL}}
\newcommand{\TT}{\mathbb T}
\newcommand{\Ga}{\mathscr G}
\newcommand{\bT}{\mathbf T}

\title{A distinction criterion for  Iwahori-spherical representations}
\author{Paul  Broussous}

\date{\today}

\begin{document}
\maketitle

\begin{abstract}
Let $G/H$ be a Galois symmetric space for an unramified quadratic extension of a locally compact field $F$,  where the group $H$ is semisimple,  simply connected, defined and split over $F$. We prove that there exists a subgroup $\Gamma = \Gamma (G/H)$ of the group of invertible elements of the Iwahori-Hecke algebra $\HIw$ of $G$ such that an Iwahori-spherical representation of $G$ is $H$-distinguished if and only if the corresponding Iwahori-Hecke module is "$\Gamma$-distinguished".
\end{abstract}

 \tableofcontents

 \bigskip

 \section{Introduction} Let $F$ be a non archimedean local field and $G/H$ be a reductive symmetric space over $F$. In particular $G$ is the group of $F$-rational points of some reductive algebraic group defined over $F$, equiped with a $F$-rational involution $\theta$, and $H$ is the subgroup of $\theta$-fixed elements.  Recall that an irreducible smooth complex representation $(\pi ,\VV )$ of $G$ is said to be {\it $H$-distinguished} if $\VV$ has a non-zero $H$-invariant linear form. These $H$-distinguished representations of $G$ are those that contribute to harmonic analysis on the homogeneous space $G/H$.  Studying these representations is the local counterpart of the global theory of periods of automorphic forms. 
\medskip

 Thanks to the theory of the Bernstein center, one knows that the category $\Sm (G)$ of smooth complex representations of $G$ splits into a direct product of indecomposable subcategories. The simplest block is that of the unramified principal series (cf. \cite{Ca}). We know since A. Borel (cf. \cite{Bo}) that this category is formed of those representations generated by their vectors fixed by an Iwahori subgroup $I$ of $G$ (we shall call these representations {\it Iwahori-spherical}). In the language of Bushnell and Kutzko's type theory, the trivial character of an Iwahori subgroup  is a type for the simplest block of $G$.  As a consequence this block is equivalent to the category of left-modules over $\HIw$, the Iwahori-Hecke algebra of $G$.
\medskip

 Let $(\pi ,\VV )$ be an irreducible Iwahori-spherical representation of $G$ and let $M$ denote the corresponding simple module over $\HIw$. We address the following question : {\it does the fact that $(\pi ,\VV )$ is $H$-distinguished correspond to a simple property of the module $M$ ? } We give an affirmative answer in the case where $G/H$ is a Galois symmetric space with some particular conditions. We indeed assume in this work that $H$ is a semisimple and simply connected group defined and split over $F$, and that $G$ is the restriction of scalars ${\rm Res}_{E/F}H$, for some unramified quadratic extension $E/F$. 
\medskip

 In the previous work \cite{Br-I}, we proved that there is a natural injective linear map $\bT$~: ${\rm Hom}_H (\VV ,\CC )\lra {\rm Hom}_\CC (M,\CC )$.  Here we determine the image of this map by proving the following result. 

\begin{theorem} (Cf. Theorem \ref{main})  There exists a subgroup $\Gamma$ of the group of invertible elements $\HIw^\times$, depending only on the Galois symmetric space $G/H$  such that for any Iwahori-spherical  irreducible representation $(\pi ,\VV )$ of $G$, the map $\bT$ induces an isomorphism 
$$
{\rm Hom}_H (\VV ,\CC )\simeq {\rm Hom}_{\Gamma} (M,\CC )\ .
$$
\end{theorem}

In other words, $\VV$ is $H$-distinguished if and only if $M$ is "$\Gamma$-distinguished". 
\smallskip

 As in \cite{Br-I} our proof is based on a geometrical model of the representation $(\pi ,\VV )$ due to Borel \cite{Bo} : $\VV$ is isomorphic to a certain space of functions defined on the chambers of the Bruhat-Tits building $X_E$ of $G$, with values in $M$. The group $\Gamma$ is given by generators of geometrical flavour. 

 Let $X_F$ be the Bruhat-Tits building of $H$, that we consider as a subset of $X_E$.  For a chamber $C$ of $X_E$, a gallery $\Ga$ in $X_E$ is called $\theta$-admissible with terminal  chamber $C$ if it is a gallery connecting a chamber $C_0$ in $X_F$ to $C$ and has minimal length among such galleries. To any $\theta$-admissible gallery $\Ga$ we associate an invertible element $e_\Ga$ of $\HIw$. We then define $\Gamma$ to be the subgroup of $\HIw^\times$ generated by the $e_{\Ga_2}^{-1}\star e_{\Ga_1}$, where $(\Ga_1 ,\Ga_2 )$ runs over the pairs of $\theta$-admissible galleries sharing the same terminal chamber. 

\medskip

 Our criterion is not effective in the sense that we do not have any effective description of the group $\Gamma$, except when $H$ has rank $1$ where in that case $\Gamma$ is trivial.  For $H ={\rm SL}(3)$, $\Gamma$ is not trivial. 
The question of describing $\Gamma$ more precisely will hopefully be the subject of future works. 
 \bigskip

 The idea of studying distinguished representations using type theory, that is using modules over Hecke algebras, was suggested to me  a few years ago by Nadir Matringe during a train journey to CIRM in Marseille. In want to thank him here, as well as Dipendra Prasad for stimulating discussions.

\section{Notation}  We try keeping the notation as in \cite{Br-I}, but for the reader's convenience we recall most of it. 
\medskip

 Let  $E/F$ an unramified quadratic extension of locally compact non-discrete non-archimedean fields. We make no restriction on the residue characteristic of $F$. We denote by $q_o$ (resp. $q$) the size of the residue field of $F$ (resp. of $E$),  so that $q=q_o^2$.

 Let $\HH$ be an algebraic group defined over $F$ that we assume connected, semisimple and simply connected. We set $\GG = {\rm Res}_{E/F}\HH$ (Weil's restriction of scalars), $H=\HH (F)$, $G=\GG (F)=\HH ( E)$. We let $X_F$ and $X_E$ denote the semisimple Bruhat-Tits buildings of $\HH$ and of $\GG$ over $F$ respectively. As polysimplicial complexes $X_F$ and $X_E$ have dimension $r={\rm rk}(\HH )$, the reductive rank of $\HH$.  Recall that $X_E$ canonically identifies with the building of $\HH$ over $E$. 

 We denote by $\theta$ the non-trivial element of ${\rm Gal}(E/F)$ and by the same letter $\theta$ the induced actions of $\theta$ on $G$ and $X_E$. Recall that there is a canonical $H\ltimes {\rm Gal}(E/F)$-equivariant embedding $X_F\lra X_E$ which is simplicial and maps chambers to chambers. We consider it as an inclusion. 
 
 We fix: 
 \smallskip
 
 -- a maximal $F$-split $\TT_o$
torus in $\HH$. We denote by $\Phi (\HH , \TT_o )$ the relative root system, $\WAff$ the associated affine Weyl group.

 -- a chamber $C_o$ in the apartment $\Ap_o$ of $X_F$ attached to $\TT_o$. We denote by $S\subset \WAff$ the set of reflexions relative to the affine hyperplans of $\Ap_o$ containing the codimension $1$ facets of $C_o$. Recall that $(\WAff ,S)$ is a Coxeter system. 
 \smallskip
 
  Observe that  $\Ap_o$ is also the apartment of $X_E$ attached to ${\TT_o^E} := {\rm Res}_{E/F}\, \TT_o$. As an apartment of $X_E$ it also has $\WAff$ as affine root system. We denote by $I := {\rm Stab}_{G} (C_o )$ the connected Iwahori subgroup of $G$ fixing $C_o$.
\medskip

 We fix the Haar measure $\mu$ on $G$ normalized by $\mu (I) =1$. We denote by $\HIw$ the convolution Hecke algebra formed of those complex functions on $G$ which are bi-$I$-invariant and have compact support, aka the {\it Iwahori-Hecke algebra} of $G$. 
 For $w\in \WAff$, we let $e_w$ denote the characteristic function of $I w I$. Recall that the $e_s$, $s\in S$, generate $\HIw$ as a $\CC$-algebra
 and that we have the quadratic relations $(e_s +1)\star (e_s -q)=0$, $s\in S$. In particular, in any $\HIw$-module, and for
 any $s\in S$, $e_s$ has its spectrum contained in $\{-1 ,q\}$.

\section{The action of $H$ on the chambers of $G$}

We denote by $\Ch_F$ (resp. by $\Ch_E$) the set of chambers in $X_F$ (resp. in $X_E$), that is the set of polysimplices of maximal dimension. If $D$ denotes a codimension $1$ facet in $\Ch_E$ (a polysimplex of dimension $r-1$), we denote by $\Ch_D$ the set of chambers in $X_E$ containing $D$. 
\smallskip

  Recall that $X_F$ is labeled by the set $S$, that is there is a map $\lambda$ : $X_F^0\lra S$, from the set of vertices of $X_E$ to $S$, such that for any chamber $C$, the restriction of $\lambda$ to the vertices of $C$ is a bijection.
  If $s$ is a vertex of $X_E$, $l(s)$ is called the {\it type} of $s$. On says that a codimension facet $D$ is {\it of type $s$}, if $\lambda (D)= S\backslash \{ s\}$. 

 One says that two chambers $C_1$ and $C_2$ are {\it adjacent}, denoted $C_1 \sim C_2$,  if $C_1\cap C_2$ has codimension $1$. We set $C_1 \sim_s C_2$, for an $s\in S$, if $C_1 \sim C_2$ and $C_1\cap C_2$ is of type $s$.

\smallskip
 
  Recall that a {\it gallery} $\Ga$ in $X_E$ is a finite sequence $(C_0 ,...,C_l )$ of chambers in $\Ch_E$ such for $i=0,...,l-1$, we have $C_i \sim C_{i+1}$. The integer $l =l(\Ga )$ is called the {\it length} of $\Ga$, and one says that $\Ga$ {\it connects $C_0$ to $C_l$}. If $C$, $C'$ are chambers in $X_E$, the {\it (combinatorial) distance} $d(C,C')$ from $C$ to $C'$ is the minimum of the $l(\Ga )$, where $\Ga$ runs over the galleries linking $C$ to $C'$. A gallery $\Ga = (C_0 ,...,C_l )$ is said to be {\it taut} if $d(C_1 ,C_2) = l(\Ga )$. If $C\in \Ch_E$, its  {\it (combinatorial) distance} $d(C,X_F)$ to $X_F$ is by definition ${\rm Min}\, \{ d(C,C') \ ; \ C'\in \Ch_F\}$. Note that if $d(C,X_F)=d$, and if $\Ga =(C_0 ,...,C_d = C)$ is a gallery with $C_0\in \Ch_F$, then for all $i=0,...,d$, $d(C_i ,X_F) = i$.

The following result, due to Fran\c cois Court\`es, is of the highest importance for us. 

\begin{theorem} \label{courtes} (F. Court\` es, \cite{BC} Proposition A.1 and its proof.) Let $D$ be a codimension $1$ facet in $X_E$. 
	
	(i) There exists an integer $d\geqslant 0$ such that for all $C\in \Ch_D$, $d(C,X_F )\in  \{ d , d+1\}$. For $k=d$, $d+1$, set $\Ch_D^k =\{ C\in \Ch_D\ ; \ d(C,X_F ) = k\}$. 
	
	(ii) The group $H$ acts transitively on $\Ch_D^d$ and $\Ch_D^{d+1}$. 
	
	(iii) The pair of cardinals  $(\vert \Ch_D^d\vert ,\vert \Ch_D^{d+1} \vert )$ is either $(q_o +1, q_o^2 -q_o )$ or $(1, q_o^2 )$.
	
	(iv) If $D\subset X_F$, we have $(\vert \Ch_D^d\vert ,\vert \Ch_D^{d+1} \vert ) =
    (q_o +1, q_o^2 -q_o )$. In general $(\vert \Ch_D^d\vert ,\vert \Ch_D^{d+1} \vert )$ depends only on the $H$-orbit of $D$. 
    
\end{theorem} 

 For a codimension $1$ facet $D$ at distance $d$ from $X_F$, we shall use the notation: $d=d(D)$,  $\vert D\vert^- =\vert \Ch_D^d$ and $\vert D\vert =\vert \Ch_D^{d+1}$.

\section{Borel's models for Iwahori-spherical representations}

Basic references for this section are \cite{Bo} and \cite{Ca}, also see \cite{Br-I}{\S}3. 
\medskip

\medskip

 We let $\Sm (G)$ denote the category of smooth complex representations of $G$ and $\Sm (G)_{I}$ denote the full subcategory of those representations  $(\pi ,\VV )$  that are generated by $\VV^{I}$, the subspace of $I$-fixed vectors.  Let $\HIw -{\rm Mod}$ denote the category of left $\HIw$-modules. Recall that the functor $\MM$ : $\Sm (G)_{I} \lra \HIw -{\rm Mod}$, $(\pi ,\VV )\mapsto \VV^{I}$ is an equivalence of categories.  We consider the pseudo-inverse $\bV$ ~: $\HH -{\rm Mod} \lra \Sm (G)_{I}$, $M\mapsto (\pi_M ,\VV_M )$ considered by Borel in \cite{Bo}{\S}4. Hence for a given representation $(r, M)$ of $\HIw$ in a $\CC$-vector space $M$, $\bV (M) = (\pi_M ,\VV_M ) $ is the corresponding representation of $G$ through the equivalence of categories.  Recall that the $G$-set $\VV_M$ is given by the tensor product:
$$
\VV_M =\HIw (G)\otimes_\HIw M,
$$
\noi where $\HIw (G)$ is the Hecke algebra\footnote{To define its convolution product, the Haar measure $\mu$  on $G$ is normalized by $\mu (I)=1$.} of $G$, considered as a right $\HIw$-module via the containment $\HIw \subset \HIw (G)$. 
\medskip

 We shall be interested in the intertwing space ${\rm Hom}_H \, (\VV_M ,\CC )$. 
Let $(r,M)$ be a representation of $\HIw$\footnote{We shall consider $M$ as a left $\HH$-module by the formula $\varphi \cdot m =r(\varphi )(m)$, $\varphi\in \HH$, $m\in M$.}. Denote by $({\tilde r}, {\tilde M})$ its contragredient, where ${\tilde M}={\rm Hom} (M,\CC )$. We define $\FF (\Ch_E , {\tilde M})_r$ to be the space of ${\tilde M}$-valued functions on $\Ch_E$ satisfying :
\begin{equation}
\label{eqn:lien}
{\tilde r}(e_s ) f (C) = \sum_{C'\sim_s C} f(C) , \ C\in \Ch_E , \ s\in S\ .
\end{equation}

 The following result is a consequence of \cite{Bo}.

\begin{proposition} (\cite{Br-I} Proposition 3.3) Let $(r,M)$ be a representation of $\HIw$. As a $\CC$-vector space  ${\rm Hom}_H \, (\VV_M,\CC )$ is isomorphic to the space $\FF (\Ch_E , {\tilde M})_r^H$ of $H$-invariant functions in $\FF (\Ch_E , {\tilde M})_r$. 
\end{proposition}

 Let $(r,M)$ be a representation of $\HIw$. Since $H$ acts transitively on $\Ch_F$, if $f\in \FF (\Ch_E , {\tilde M})_r^H$ the value $f(C)$ does not depend on $C\in \Ch_F$. We define a linear map $\bT$~: $\FF (\Ch_E , {\tilde M})_r^H \lra {\tilde M}$ by $\bT (f)=f(C_o )$, where $C_o$ is our fixed chamber in $\Ch_F$. We recall the fundamental fact:

\begin{theorem} (\cite{Br-I} Proposition 4.2). The map $\bT $~: $\FF (\Ch_E , {\tilde M})_r^H \lra {\tilde M}$ is injective. 
\end{theorem} 

 Let us note that in Proposition 4.2 of \cite{Br-I}, $M$ is assumed irreducible. However this assumption is not used in its proof. 

    In the next sections we shall determine the image of $\bT$.

\section{Image of $\bT$: a necessary condition}

 We fix a representation $(r, M)$ of $\HIw$.  Let $f$ be a $H$-invariant fonction $\Ch_E \lra {\tilde M}$.  Assume that $f$ satisfies condition \ref{eqn:lien}. Let $D$ be a codimension $1$ facet of $X_E$, set $d=d(D)$ and let $s\in S$ be the type of $D$. By Theorem \ref{courtes}, $f$ is constant on $\Ch_D^{d}$ and $\Ch_D^{d+1}$. We let $f^+ (D)$ (resp. $f^- (D)$) denote the constant value of $f$ on $\Ch_D^d$ (resp. on $\Ch_D^{d+1}$). 
 If $C$ is a chamber in $\Ch_D$, taking the $H$-invariance of $f$ into account, Relation \ref{eqn:lien} takes the following simple form:
\smallskip
 
  {\it Case no 1}. If $C\in \Ch_D^d$, then Relation \ref{eqn:lien} writes
  
\begin{equation} \label{eqn:case1}	
f^+ (D) = \frac{1}{\vert D\vert^+}  {\tilde r} (e_s -(\vert D\vert^- -1)e_1)\, f^{-}(D)
\end{equation}

 {\it Case no 2}. If 
 $C\in \Ch_D^{d+1}$, then Relation \ref{eqn:lien} writes
 
 \begin{equation} \label{eqn:case2}
 	f^- (D) = \frac{1}{\vert D\vert^-} {\tilde r} (e_s -(\vert D\vert^+ -1)e_1)\, f^{+}(D)
 \end{equation}

 So a $H$-invariant function $f$~: $\Ch_E\lra {\tilde M}$ satisfies Relation \ref{eqn:lien} for all $C\in \Ch_E$ and all $s\in S$ if, and only if, $f$ satisfies Relations \ref{eqn:case1} and \ref{eqn:case2} for all codimension $1$ facet $D$. In fact  Relations \ref{eqn:case1} and \ref{eqn:case2} are equivalent.

 \begin{proposition}  \label{liensimple} A $H$-invariant function satifies \ref{eqn:lien} for all $s$ and $C$ if and only if it satisfies \ref{eqn:case1}. 
 \end{proposition}

\noi {\it Proof}. The proposition is a straightforward consequence of the following lemma. Its proof, left to the reader, is a simple calculation based on the quadratic relations in $\HIw$ (note that for all codimension $1$ facet $D$ in $X_E$, we have $\vert D\vert^- +\vert D\vert^+ =q+1$). 

\begin{lemma} Let $a$, $b$ be non-zero complex numbers such that $a+b=q+1$. Then for all $s\in S$, $\frac{1}{a} \, (e_s -(b-1)e_1 )\in \HIw$ is invertible with inverse $\frac{1}{b}\, (e_s -(a-1)e_1 )$. 
\end{lemma}

 To any   codimension $1$ facet $D$ in $X_E$ we attach the following invertible element of $\HIw$:
 
 \begin{equation} \label{eqn:eD}
 	e_D := \frac{1}{\vert D\vert^+}\, (e_s -(\vert D\vert^- - 1)e_1 )\in \HIw^\times
\end{equation}
 	
 \noi where $s$ is the type of $D$.  Since the action of $H$ on $X_E$  preserves	the labelling, we have $e_{h.D}=e_D$,  for all $h\in H$.
 \medskip
 
  A gallery $\Ga = (C_0 ,...,C_d)$ in $X_E$ is called $\theta$-admissible if $C_0 \in \Ch_F$ and $d=d ( C_d ,X_F )$. If $\Ga$ is such a gallery, we set $D_i = C_i \cap C_{i+1}$, $i=0,...,d-1$. We call $C_0$ (resp. $C_d$) the {\it initial} chamber (resp. the  {\it terminal} chamber) of $\Ga$. We shall use the notation:

 	\begin{equation}\label{eqn:eG}
 	 e_{\Ga} := e_{D_{d-1}} \star e_{D_{d-2}} \star \cdots \star e_{D_0}  \in \HIw^\times\ .
 	\end{equation}
 
  If $d=0$, by convention we set $e_{\Ga}=e_1$. Let us observe that $e_{h\Ga} = e_{\Ga}$ for any $\theta$-admissible gallery $\Ga$ and all $h\in H$.

  \begin{proposition} \label{necessaire} Let ${\tilde m}\in {\tilde M}$. If $\tilde m$ lies in the image of $\bT$, then for all pairs of $\theta$-admissible galleries $\Ga_1$ and $\Ga_2$ with the same terminal chamber, we have $e_{\Ga_1} \cdot {\tilde m} = e_{\Ga_2} \cdot {\tilde m}$. 
  \end{proposition} 

\noi {\it Proof}. Write ${\tilde m}= \bT (f)$, with $f\in \FF (\Ch_E , {\tilde M})_r^H$. By definition ${\tilde m}=f(C)$, for all $C\in \Ch_F$. Let $C\in \Ch_E$ and let $\Ga$ any $\theta$-admissible gallery with terminal chamber $C$. By Relation \ref{eqn:case1}, we have $f(C) = e_{\Ga}\cdot {\tilde m}$. Hence $e_{\Ga}\cdot {\tilde m}$ depends only on the terminal chamber of $\Ga$; the proposition follows.

\section{A criterion}

 We start by proving a converse to Proposition \ref{necessaire}.

\begin{proposition} \label{reciproque}  Let $(r,M)$ be a representation of $\HIw$ and ${\tilde m}\in {\tilde M}$. Assume that for all pairs of $\theta$-admissible galleries $\Ga_1$, $\Ga_2$ with the same terminal chamber, we have $e_{\Ga_1}\cdot {\tilde m} = e_{\Ga_2}\cdot {\tilde m}$. Then $\tilde m$ lies in the image of $\bT$.
\end{proposition}

 \noi {\it Proof}. Let ${\tilde m}\in {\tilde M}$ satisfying the assumption of the proposition. Define a function $f$~: $\Ch_E\lra {\tilde M}$ by $f(C) = e_{\Ga}\cdot {\tilde m}$, for all $C\in \Ch_E$, where $\Ga$ is any $\theta$-admissible gallery with terminal chamber $C$. Thanks to our assumption, $f(C)$ does not depend on the choice of $\Ga$. 

 Let $C\in \Ch_E$ and $h\in H$. Let $\Ga$ be an $\theta$-admissible gallery with terminal chamber $C$. Then $h\cdot \Ga$ is a $\theta$-admissible gallery with terminal chamber $h.C$. Since $e_{h.\Ga}=e_{\Ga}$, we have that $f(h.C)=f(C)$ and $f$ is $H$-invariant. 

We now prove that $f\in \FF (\Ch_E ,{\tilde M})_r^H$. For this we prove that $f$ satisfies Relation \ref{eqn:lien}, or equivalently satisfies Relation \ref{eqn:case1}.  Let $D$ be any codimension $1$ facet and set $d=d(D)$; we may assume that $d\geqslant 1$.  Let $C\in\Ch_D^{d+1}$. Choose an $\theta$-admissible gallery $\Ga = (C_0 ,C_1 ,..., C_d )$ with $C_d =C$. Then $C_{d-1}\in \Ch_D^d$. Set $D_i := C_i \cap C_{i+1}$, $i=0,...,d$, so that $D_{d-1}=D$.  Let $\Ga ' = (C_0 ,...,C_{d-1})$; this is an $\theta$-admissible gallery with terminal chamber $C_{d-1}$.  We have $e_{\Ga}=e_{D}\star e_{\Ga '}$. By definition of $f$, we have:
\smallskip

 (a) $f^+ (D)= f(C) = e_{\Ga }\cdot {\tilde m} = e_{D} \cdot (e_{\Ga '} \cdot {\tilde m})$, 

 (b) $f^- (D) = f(C_{d-1})=e_{\Ga '}\cdot {\tilde m}$. 
\smallskip

 It follows that $f^+ (D)=e_{D}\cdot f^{-} (D)$, and that $f$ satisfies Relation \ref{eqn:case1}. 
\smallskip

 Now by construction we have $\bT (f)={\tilde m}$, as required.

\medskip

 Let $\Gamma$ be the subgroup of $\HIw^\times$ generated by the $e_{\Ga_2}^{-1}\star e_{\Ga_1}$, where $(\Ga_1 ,\Ga_2 )$ runs over the pairs of $\theta$-admissible galleries in $\Ch_E$ sharing the same terminal chamber. This subgroup $\Gamma$ depends only on the symmetric pair $(G,H)$. 
\medskip

Clearly Proposition \ref{reciproque} may be rephrased as follows. 

\begin{corollary}  Let $(r,M)$ be a representation of $\HIw$. Then the image of 
$\bT$ is ${\tilde M}^\Gamma$.
\end{corollary}

We deduce the following criterion for the $H$-distinction of an Iwahori-spherical  representation of $G$  in terms of its module over $\HIw$.

\begin{theorem}\label{main}  Let $(\pi ,\VV )$ be an Iwahori-spherical representation of $G$ with Iwahori-Hecke algebra module $M$. Then we have an isomorphism of $\CC$-vector spaces:
$$
{\rm Hom}_H (\VV ,\CC ) \simeq {\tilde M}^\Gamma = {\rm Hom}_\Gamma (M,\CC )
$$
where ${\tilde M}^\Gamma$ is the space of $\Gamma$-fixed elements in the $\HIw$-module $\tilde M$.

In other words $(\pi ,\VV )$ is $H$-distinguished if, and only if, ${\rm Hom}_\Gamma (M,\CC ) \not= 0$, and in this case the multiplicity is ${\rm dim}\, {\rm Hom}_\Gamma (M,\CC )$. 
\end{theorem}

Let us give an example. Assume that the $F$-rank of $\HH$ is $1$, or equivalently that $X_F$ and $X_E$ are trees. Since $\HH$ is semisimple and simply connected, we necessarily have $H={\rm SL}(2,F)$ and $G={\rm SL}(2,E)$. Let $C$ be a chamber in $\Ch_E$ that does lie in $\Ch_F$.  Write $C=\{ a,b\}$, where the vertices $a$ and $b$ are such that $d(b,X_F ) > d(a,X_F)$.  Then the enclos of $\ds C\cup C^\theta$ is a geodesic segment in $X_E$. It is contained in any ${\rm Gal}(E/F)$-stable appartment containing $C$.  Let $x_C$ denote the intersection of that geodesic segment with $X_F$. The vertex $x_C$ is the projection of any point of $C$ onto the convex closed subset $X_F\subset X_E$. Any $\theta$-admissible gallery with terminal chamber $C$ has the form $\Ga = (C_0, C_1 ,...,C_d )$, where : 
\smallskip

 (a) $C_0$ is one of $q_o +1$ chambers of $X_F$ containing $x_C$,

 (b) $C_1 ,...,C_d$ are the chambers of $X_E$ contained in the geodesic segment $[x_C ,b]$.
\smallskip

 It follows that the invertible element $e_\Ga$ of $\HIw$ only depends on $C$.  So the group $\Gamma$ is trivial in this case:

\begin{corollary} Any Iwahori-spherical irreducible representation $(\pi ,\VV )$ of ${\rm SL}(2,E)$ is ${\rm SL}(2,F)$-distinguished with multiplicity ${\rm dim}\, \VV^{I}$.
\end{corollary}

 It is not difficult to construct non-trivial elements of $\Gamma$ when $H={\rm SL}(3)$. In general the group $\Gamma$ seems difficult to determine.  To have more information on $\Gamma$ we need a better understanding of the structure of $\theta$-admissible galleries in $X_E$.

Paul Broussous
\smallskip

paul.broussous{@}math.univ-poitiers.fr
\medskip

 Laboratoire de Math\'ematiques et Applications, UMR 7348 du CNRS

Site du Futuroscope -- T\'el\'eport 2
 
11,  Boulevard Marie et Pierre Curie

B\^atiment H3 --  TSA 61125

86073 POITIERS CEDEX 

 France


\begin{thebibliography}{99}



\bibitem{Bo} A. Borel, {\it Admissible representations of a 
semi-simple group over a local field with vectors fixed under 
an Iwahori subgroup},   Invent. Math.  35  (1976), 233--259.






\bibitem{BC} P. Broussous and F. Court\`es,  {\it Distinction of the
  Steinberg representation,  with an appendix by Fran\c cois Court\`es.},  Int. Math. Res. Not. IMRN 2014, no. 11, 3140–3157.

\bibitem{Br-I} P. Broussous, {\it On the distinction of Iwahori-spherical discrete series representations},   Journal of the Institute of Mathematics of Jussieu. Published online 2024:1-15. doi:10.1017/S1474748024000185 


\bibitem{Ca} W. Casselman, {\it The unramified principal series of p-adic groups. I. The spherical function}, Compositio Math. 40 (1980), No. 3, 387--406. 

\bibitem{Cou} F. Court\`es, {\it Distinction of the Steinberg representation II: an equality of characters},  Forum Math. 27 (2015), no. 6, 346--3475.














\end{thebibliography}
\end{document}